\documentclass[english]{ourlematema}

\usepackage[T1]{fontenc}
\usepackage[english]{babel}
\usepackage{todonotes}

\usepackage{mathtools}
\usepackage{amsmath}                                                            
\usepackage{amssymb}
\usepackage{amsthm}
\usepackage{amscd}
\usepackage{algorithm, algorithmic}
\usepackage{listings}

\usepackage{caption}
\usepackage{subcaption}

\usepackage{extarrows}

\usepackage[colorlinks, linktocpage, citecolor = red, linkcolor = blue]{hyperref}
\usepackage{color}

\usepackage[shortlabels]{enumitem}

\usepackage[capitalise,noabbrev]{cleveref} 

\newtheorem{theorem}{Theorem}[section]
\newtheorem{theorem*}{Theorem*}[section]

\theoremstyle{definition}

\newtheorem{definition}[theorem]{Definition}

\newtheorem{example}[theorem]{Example}

\newcommand{\C}{\mathbb{C}}

\newcommand{\CC}{\mathbb{C}}

\newcommand{\PP}{\mathbb{P}}

\newcommand{\PC}{\PP_{\CC}}

\usepackage{tikz-cd}
\tikzcdset{
arrow style       = tikz,
cramped,
diagrams          = {>=stealth},
row sep/normal    = 2.2em,
column sep/normal = 2.7em,
}

\renewenvironment{tikzcd}{\begin{oldtikzcd}}{\end{oldtikzcd}\vspace{-3ex}}

\newcommand*{\classicsets}[1]{\mathbb{#1}} 

\newcommand*{\comp}{\classicsets{C}}   

\newcommand*{\rightstealtharrow}{%
  \mathrel{\begin{oldtikzcd}[ampersand replacement=\&, column sep=1.15em]\arrow[r]\&{}\end{oldtikzcd}}}

\newcommand*{\myapplication}[1]{\mspace{-8mu} #1 \mspace{-8mu}} 
\renewcommand*{\to}{\myapplication{\rightstealtharrow}}

\DeclareMathOperator{\pgl}{PGL}

\title{96120: The degree of the linear orbit \\ of a cubic surface}
\titlemark{Degree of the linear orbit closure }

\titlenote{We thank Elisa Cazzador, Corey Harris, Kristian Ranestad, Anna Seigal, Bjorn Skauli, and Bernd Sturmfels for helpful discussions, and Bernd Sturmfels for suggesting the problem.
  Research on this project was carried out in part while the authors were based at the Max Planck Institute for Mathematics in the Sciences (MPI-MiS) in Leipzig, Germany.
  LBM was partially supported by the Mineco grant MTM2016-75980-P.
  ST was supported by the Deutsche Forschungsgemeinschaft (German Research Foundation) Graduiertenkolleg {\em Facets of Complexity} (GRK~2434).
This material is based upon work supported by the National Science Foundation Graduate Research Fellowship Program under Grant No. DGE 1752814. Any opinions, findings, and conclusions or recommendations expressed in this material are those of the authors and do not necessarily reflect the views of the National Science Foundation.}


\author{Laura Brustenga i Moncus\'i}
\address{%
  Departament de Matem\'atiques\\
  Universitat Aut\'onoma de Barcelona\\
  08193 Bellaterra, Spain\\
  \email{\href{mailto:brust@mat.uab.cat}{brust@mat.uab.cat}}
  }

\author{Sascha Timme}
\address{%
  Institut f\"ur Mathematik\\
  Technische Universit\"at Berlin\\
  Stra{\ss}e des 17. Juni 136, 10623 Berlin\\
  \email{\href{mailto:timme@math.tu-berlin.de}{timme@math.tu-berlin.de}}
  }
\authorline
\author{Madeleine Weinstein}
\address{%
  Department of Mathematics\\
  University of California Berkeley\\
  970 Evans Hall, Berkeley, CA 94720\\
  \email{\href{mailto:madeleine\_weinstein@berkeley.edu}{ madeleine\_weinstein@berkeley.edu}}
  }

\shortauthormark

\date{2019/09/14}

\begin{document}

\maketitle

\begin{abstract}
The projective linear group \(\pgl(\comp,4)\) acts on cubic surfaces, considered as points of $\PC^{19}$. We compute the degree of the $15$-dimensional projective variety given by the Zariski closure of the orbit of a general cubic surface. The result, 96120, is obtained using methods from numerical algebraic geometry.
\end{abstract}

\maketitle

\section{Introduction}
Automorphism groups of varieties and group actions on varieties are of much interest to researchers of algebraic geometry, arithmetic, and representation theory \cite{Aluffi:Faber:1993, matserror,Mats,cayleycubic}.
Here we study the action of the projective linear group  \(\pgl(\comp,4)\) on cubic surfaces parameterized by points in $\PC^{19}$.  In particular, we compute the degree of the $15$-dimensional projective variety in $\PC^{19}$ defined by the Zariski closure of the orbit of a general cubic surface under this action. 
This degree is also meaningful in enumerative geometry: It is the number of translates of a cubic surface that pass through 15 points in general position. This formulation provides an alternate method for obtaining the degree.

Aluffi and Faber considered the analogous problem for plane curves of arbitrary degree, first the smooth case in \cite{Aluffi:Faber:1993} and then the general case in \cite{Aluffi:Faber:2000}.
They obtained a closed formula for the degree of the orbit closure of a plane curve under the action of \(\pgl(\comp,3)\).
This was a significant undertaking, involving long and detailed calculations in intersection rings using advanced techniques from intersection theory.

Instead of adopting the techniques developed by Aluffi and Faber, we use tools from \emph{numerical algebraic geometry} \cite{Hauenstein:Sommese:2017,Sommese:Wampler:2005}.
The general idea is as follows.
We fix a cubic surface $f$ and $15$ points in general position in $\PC^3$. 
The condition that a translate of $f$ passes through these 15 points results in a polynomial system for which we compute all isolated numerical solutions by homotopy continuation and monodromy methods using the software {\tt HomotopyContinuation.jl} \cite{HomotopyContinuation.jl}.
The concept of an \emph{approximate zero} \cite{BCSS:1998} makes precise the definition of a numerical solution.
We use Smale's $\alpha$-theory and the software {\tt alphaCertified} \cite{Hauenstein:Sottile:2012} to certify that the obtained numerical solutions indeed satisfy the system of polynomial equations.
Finally, we use a \emph{trace test} \cite{Leykin:Rodriguez:Sottile:2018} to check that no solution is missing.
With these techniques, we conclude that the number of numerical solutions we obtain, 96120, is in fact the degree of the orbit closure. This result is a ``numerical theorem'' rather a theorem in the classical sense.
\medskip

Our presentation is organized as follows.
In Section \ref{sec:orbits}, we introduce the linear orbit problem in detail and derive the polynomial systems used in our computations.
In Section \ref{sec:numerical-algebraic-geometry}, we discuss the techniques used from numerical algebraic geometry and in Section \ref{sec:results} the computations performed to arrive at the result.

\section{Linear Orbits and Polynomial Systems}\label{sec:orbits}
A cubic surface in $\PC^3$ is defined by a cubic homogeneous polynomial in $4$ variables with complex coefficients.
The parameter space for cubic surfaces is $\PC^{19}$ and we fix coordinates \((c_0:\dots:c_{19})\in\PC^{19}\).

The projective space $\PC^{15}$ of homogeneous $4 \times 4$ matrices $A=(a_{ij})_{1 \le i,j \le 4}$ is a compactification of the projective general linear group
$$
\pgl(\CC, 4)=\{ A \in \PC^{15} \,|\, \det A \neq 0 \}\subseteq \PC^{15}\, .
$$
The group \(\pgl(\comp,4)\) acts on a cubic surface $f\in \PC^{19}$, with \(\varphi\in \pgl(\comp,4)\) sending $f$ to the cubic surface $\varphi \cdot f$ defined by the equation
$$f(\varphi(x,y,z,w)) = 0 \,.$$
This corresponds to a linear change of the coordinates \(x,y,z,w\). We say that $\varphi \cdot f$ is the \emph{translate} of $f$ by $\varphi$.
Then $\pgl(\comp,4) \cdot f$ is the orbit of $f$ in $\PC^{19}$ and its Zariski closure $\Omega_f := \overline{\pgl(\CC, 4) \cdot f}$ is a 15-dimensional projective variety.

\begin{example}
  To illustrate this idea, we consider the action of $\pgl(\mathbb{C},2)$ on pairs of points defined by homogeneous polynomials
  $$f(x,y)=b_0x^2+b_1xy+b_2y^2\,.$$
  The parameter space for pairs of points is $\PC^{2}$, that is $f=(b_0:b_1:b_2) \in \PC^2$. 
  Let
  $$\varphi= \begin{pmatrix}
  a_{11} & a_{12} \\
  a_{21} & a_{22}
  \end{pmatrix}. $$
  Then
  \begin{align*}
  f(\varphi(x,y))=& b_1(a_{11}x+a_{12}y)^2+b_2(a_{11}x+a_{12}y)(a_{21}x+a_{22}y)+b_3(a_{21}x+a_{22}y)^2 \\ =&  (b_1a_{11}^2+ b_2a_{11}a_{21}+b_3a_{21}^2)x^2 \;+ \\ & (2b_1a_{11}a_{12} + b_2(a_{11}a_{22}+a_{12}a_{21})+2b_3a_{21}a_{22})xy \;+ \\ & (b_1a_{12}^2+ b_2a_{12}a_{22}+b_3a_{22}^2)y^2.
  \end{align*}
  and thus
  \begin{align*}
    \varphi \cdot f = ( &b_1a_{11}^2+ b_2a_{11}a_{21}+b_3a_{21}^2 : \\
  &2b_1a_{11}a_{12} + b_2(a_{11}a_{22}+a_{12}a_{21})+2b_3a_{21}a_{22}: \\
 & b_1a_{12}^2+ b_2a_{12}a_{22}+b_3a_{22}^2)\in\PC^{2} \,.
  \end{align*}
\end{example}

To compute the degree of the orbit closure of a general cubic surface under the action of \(\pgl(\comp,4)\), we construct as follows polynomial systems whose number of isolated regular solutions correspond to the desired degree.

Fix a general cubic surface $f \in \PC^{19}$ and a general linear subspace \(L\subseteq \PC^{19}\) of dimension 4, the codimension of \(\Omega_f\).
Consider the rational map
$$
\Theta_f \,:\, \PC^{15}\to \PC^{19}
$$ sending a $4\times 4$ matrix $\varphi$ to $\varphi\cdot f$.
By definition, the image of $\Theta_f$ is $\Omega_f$. 
By \cite[Theorem 5]{Mats}, a generic hypersurface of degree at least three in at least four variables has a trivial stabilizer (we note that in \cite[Propostion 7.5]{matserror} it is stated that argument in \cite{Mats} has an error but that it does not affect the correctness of the statement).
Hence, the map $\Theta_f$ is one-to-one, so the degrees of the zero-dimensional varieties $\Omega_f  \cap L$ and $\Theta_f^{-1}(\Omega_f \cap L)=\Theta_f^{-1}(L)$ are equal.

Note that $\Theta_f^{-1}(L)$ includes non-invertible matrices whose kernel does not contain $f$.
But since we assume $L \subseteq \PC^{19}$ to be general, \(\Theta_f^{-1}(L)\) will not intersect the codimension $1$ subvariety of $\PC^{15}$ of matrices with determinant equal to $0$. 

It follows 
that the degree of the orbit closure is the number of regular isolated solutions of the polynomial system
\begin{equation}\label{eq:formulation_subspace}
  \tilde{L} \; \varphi \cdot f = 0
\end{equation}
in the entries of \(\varphi\in \PC^{15}\), where $\tilde{L} \in \C^{15 \times 20}$ is a matrix representing the general linear subspace \(L\subseteq \PC^{19}\) of dimension 4.

\medskip

The degree of \(\Omega_f\) can be thought of in enumerative terms as the number of translates of $f$ that pass through 15 points $p_1,\ldots,p_{15} \in \PC^3$ in general position.
Consider the translated cubic surface \(\varphi\cdot f\).
Note that \(\varphi\cdot f\) passes through a point $p \in \PC^3$ if and only if $f(\varphi(p)) = 0\,.$
Therefore we obtain the polynomial system
\begin{equation}\label{eq:formulation_enumerative}
  f(\varphi(p_i)) = 0 \;, \; i=1,\ldots,15
\end{equation}
in the entries of \(\PC^{15}\).
By Bertini's theorem, we may assume that the hypersurfaces $f(\varphi(p_i))=0$ intersect transversally.
Hence, the degree of $\Omega_f$ is equal to the number of matrices satisfying \eqref{eq:formulation_enumerative}.

Formulations \eqref{eq:formulation_subspace} and \eqref{eq:formulation_enumerative} both result in a system of $15$ homogeneous cubic polynomials in the 16 unknowns \((a_{ij})_{1\le i,j\le 4}\), but they have different computational advantages.
To perform numerical homotopy continuation, it is beneficial to pass to an affine chart of projective space.
This can be done in formulation \eqref{eq:formulation_subspace} by fixing a coordinate, say adding the polynomial $a_{11}-1=0$.
But this introduces artificial solutions.
For example, for every solution $\phi \in \CC^{16}$, we have that
$e^{i\frac23 \pi}\phi$ and $e^{i\frac43 \pi}\phi$ are also solutions.
The formulation \eqref{eq:formulation_enumerative} does not produce these undesired artificial solutions.
However, the formulation \eqref{eq:formulation_subspace} is better suited for applying the trace test than \eqref{eq:formulation_enumerative}.
The reason is given in the following section. 

\section{Numerical Algebraic Geometry}
\label{sec:numerical-algebraic-geometry}

Numerical algebraic geometry concerns numerical computations of objects describing algebraic sets defined over subfields of the complex numbers.
The most basic of these objects are the \emph{solution sets}, a data structure for representing solutions to polynomial systems. 
The term ``numerical'' refers to computations which are potentially inexact (e.g., floating-point arithmetic).
However, this does not necessarily mean that the results obtained are unreliable.
The certification of solutions plays an important role in the field.
For a more in-depth definition and a brief history of numerical algebraic geometry see \cite{Hauenstein:Sommese:2017}.
A comprehensive introduction to the subject is available in \cite{Sommese:Wampler:2005}.
\medskip

We now introduce tools from numerical algebraic geometry needed to compute and certify the degree of the orbit closure.
We fix a system of polynomials \(F=(F_1,\dots,F_m)\) in $n$ variables and assume that it has $l$ isolated solutions $p_1,\dots,p_l \in \CC^n$. 

\paragraph{Homotopy Continuation.} Numerical homotopy continuation \cite[Section 8.4.1]{Sommese:Wampler:2005} is a fundamental method that underlies most of numerical algebraic geometry. 
The general idea is as follows.
Suppose we want to compute the isolated solutions of $F$.
We build a homotopy $H(x,t): \CC^n \times \C \rightarrow \C^m$ which deforms a system of polynomials $G(x)=H(x,0)$ whose isolated solutions are known or easily computable into the system $F(x)=H(x,1)$. A well-defined homotopy requires that $G$ has at least as many isolated solutions as \(F\) so that we are able to compute \emph{all} isolated solutions of $F$.
Given a solution $x_0$ of $G$, there is a solution path $x(t): \CC \to \CC^n$, which is a curve implicitly defined by the conditions $x(0) = x_0$ and $H(x(t),t) = 0$ for $t \in U \subseteq \C$ where $U$ is the flat locus of the projection $\CC^n\times \C\to \C$ restricted to $H=0$, which is dense in $\C$ by generic flatness. In particular, a well-defined homotopy requires $0 \in U$.
The solution path is usually tracked using a predictor-corrector scheme.
As $t$ approaches $1$ the solution path either diverges or converges to a solution of $F$. 

A standard homotopy is the \emph{total degree homotopy}. B\'ezout's theorem gives \(N=\prod_{i=1}^m\deg(F_i)\) as an upper bound for the  the number of isolated solutions of $F$. A total degree homotopy uses a start system $G$ with $N$ isolated solutions and the homotopy $H(x,t)=(1-t)G(x)+tF(x)$. As the B\'ezout bound may be very high, for large computations the total degree homotopy is impractical and other methods are necessary.
\medskip

\paragraph{Monodromy method.}
Monodromy (see~\cite{Monodromy:18,delCampo:Rodriguez:2017}) is an alternative method for finding isolated solutions to parameterized polynomial systems which is advantageous if the number of solutions is substantially lower than the B\'ezout bound. 
Embed our polynomial system \(F\) in a family of polynomial systems $\mathcal{F}_Q$, parameterized by a connected open set $Q \subseteq \mathbb{C}^k$.
Let $l$ be the number of solutions of $F_q\in \mathcal{F}_Q$ for $q\in U$, where $U\subseteq Q$ is the flat locus of the family $\mathcal{F}_Q$.

Consider the incidence variety
$$
   Y\,:=\,\bigl\{(x,q) \in \mathbb{C}^n \times Q \; | \; F_q(x) = 0\bigr\}\,.
$$
Let $\pi$ be the projection from $\mathbb C^n\times Q$ onto the second argument restricted to \(Y\).
For every $q \in Q\backslash \Delta$, the fiber $Y_q=\pi^{-1}(q)$ has exactly $l$ points.
Given a loop \(O\) in \(U\) based at \(q\), the preimage \(\pi^{-1}(O)\) is a union of paths starting and ending at (possibly different) points of \(Y_q\).
So, giving a direction to the loop \(O\), we may associate to each point \(y\) of \(Y_q\) the endpoint of the path starting at \(y\).
This defines an action, the \emph{monodromy action}, of the fundamental group of \(U\) on the fiber \(Y_q\), 
which in turn defines a map from the fundamental group of \(U\) to the symmetric group \(S_l\).
The \emph{monodromy group} of our family at \(q\) is the image of such a map.
This action is transitive if and only if $Y$ is irreducible, which we assume.  

Fix \(q_0\in U\) such that $F=F_{q_0}\in\mathcal{F}_Q$.
Suppose a \emph{start pair} $(x_0,q_0)$ is given, that is, $x_0$ is a solution to the instance $F_{q_0}$. 
The start solution $x_0$ is numerically tracked along a directed loop in $Q\backslash \Delta$, yielding a solution $p_0'$ at the end.
If $p_0\not= p_0'$, then $p_0'$ is tracked along the \emph{same} loop, possibly yielding again a new solution. 
Then, all solutions are tracked along a \emph{new} loop, and the process is repeated until some stopping criterion is fulfilled. 

We note that this method requires us to know one solution of our polynomial system to use as a start pair. Various strategies exist to find such a solution. We will describe one strategy in Section \ref{sec:results}. 

\paragraph{Certifying solutions.} 
The above methods yield numerical approximations of solutions of our polynomial system $F$. How can we certify that the obtained approximations correspond to actual solutions of $F$ and that they are all distinct? For systems $F$ with an equal number $n$ of polynomials and variables, Smale introduced the notion of an \emph{approximate zero}, the \(\alpha\)-number and the \(\alpha\)-theorem, see \cite{Smale:86}. In short, an approximate zero of $F$ is any point $p\in \mathbb{C}^n$ such that Newton's method, when applied to $p$, converges quadratically towards a zero of $F$. This means that the number of correct significant digits roughly doubles with each iteration of Newton's method. 

\begin{definition}[Approximate zero]\label{def_approx_zero} \rm
Let
$J_F$ be the $n \times n$ Jacobian matrix of $F$.
A point $p \in \mathbb{C}^n$ is an {\em approximate zero} of $F\,$ if there exists a zero $\zeta\in \mathbb{C}^n$ of~$F$ such that the sequence of Newton iterates
$$z_{k+1} = z_k - J_F(z_k)^{-1}F(z_k)
$$ starting at $z_0=p$ satisfies for all $k\ge1$ that
$$ \Vert z_{k+1} - \zeta\Vert \,\leq \, \frac{1}{2}\Vert z_k-\zeta\Vert^2.$$
If this holds, then we call $\zeta$ the {\em associated zero} of~$p$.
Here~$\Vert {x}\Vert$
is the standard Euclidean norm in~$\mathbb{C}^n$, and the zero $\zeta$ is assumed to be nonsingular (that is, ${\rm det}(J_F(\zeta)) \not=0$ since \(F\)).
\end{definition}

To check whether a point $p\in\C^n$ is an approximate zero of $F$ from Definition \ref{def_approx_zero} requires infinitely many steps, one for each iteration of the Newton method.
Nevertheless, when \(p\) is close enough to its associated zero, it is possible to certify that $p$ is an approximate zero with only finitely many computations, as we now see. Smale's \(\alpha\)-theorem (see~\cite[Theorem 4 in Chapter 8]{BCSS:1998}) is an essential ingredient.
The theorem uses the $\gamma$- and $\alpha$-numbers
\begin{align*}
\gamma(F,x) &\,\,=\,\, \sup_{k\geq 2}\big\Vert \frac{1}{k!}\,J_F(x)^{-1} D^kF(x)\big\Vert^\frac{1}{k-1}\ \text{ and}\\
\alpha(F,x)  &\,\,=\,\,  \Vert J_F(x)^{-1}F(x)\Vert  \cdot \gamma(F,x)\,,
\end{align*}
where $D^kF$ is the tensor of order-$k$ derivatives of \(F\) and the tensor $J_F^{-1}D^kF$ is understood as a multilinear map $A:(\mathbb{C}^n)^k\to \mathbb{C}^n$ with norm $\Vert A\Vert := \max_{\Vert v \Vert = 1} \Vert A(v,\ldots,v)\Vert$.
\begin{theorem}[Smale's $\alpha$-theorem]\label{alpha_theorem}
If $\alpha(F,x)<0.03$, then $x$ is an approximate zero of $F$. Furthermore, if $y\in\mathbb{C}^n$ is any point with $\Vert y-x\Vert$ less than $(20\,\gamma(F,x))^{-1}$, then $y$ is also an approximate zero of $F$ with the same
associated zero $\zeta$ as $x$.
\end{theorem}

Smale's $\alpha$-theorem is in fact more general than is stated above. The numbers 0.03 and 20 can be replaced by any pair of positive numbers satisfying certain constraints.

To avoid the computation of the $\gamma$-number
Shub and Smale \cite{SS} derived an upper bound for $\gamma(F,x)$ which can be computed exactly and efficiently.
Hence, one can decide algorithmically whether $x$ is an approximate zero using only the data of the point~$x$ itself and \(F\). 
Hauenstein and Sottile \cite{Hauenstein:Sottile:2012} implemented these ideas
in~an algorithm, called {\tt alphaCertified}, which decides both whether a point $x\in \mathbb{C}^n$ is an approximate zero and whether two approximate zeros have distinct associated zeros.

\paragraph{Trace test} 

The certification process explained above establishes a \emph{lower} bound for the number of isolated solutions of $F$.  The trace test can be used for polynomial systems satisfying certain conditions to show that \emph{all} solutions have been found. See \cite{Leykin:Rodriguez:Sottile:2018} for a more detailed explanation. 

We first establish definitions of concepts used in the trace test. 
A \emph{pencil of linear spaces} is a family $M_t$ for $t \in \C$ of linear spaces that depends affinely on the parameter $t$. Each $M_t$ is the span of a linear space $L$ and a point $t$ on a line $l$ that is disjoint from $L$.
Suppose that $W \subset~\C^n$ is an irreducible variety of dimension $m$ and that $M_t$ for $t \in \C$ is a general pencil of linear subspaces of codimension $m$ such that $W$ and $M_0$ intersect transversally. 
Consider a fixed subset $W' \subseteq W \cap  M_0$ and denote by $W^{'}_t \subseteq W \cap M_t$ the points obtained by tracking $W'$ along the pencil.
Denote by $w(t)$ the sum of the points of $W'_t$. If $W'_t = W \cap M_t$ then $w(t)$ is the \emph{trace} of $W \cap M_t$. A $\C$-valued function $w$ is called an \emph{affine linear function} of $t$ if there exist $a, b \in \CC$ such that $w(t) = a + bt$. A $\C^n$-valued function $w$ is called an affine linear function of $t$ if for a nonconstant path $\gamma: [0,1] \to \mathbb{C}$ with $\gamma(0)=0$, we have that $w(\gamma(s))$ is an affine linear function of $\gamma(s)$. The trace is an affine linear function of $t$ \cite[Prop. 3]{Leykin:Rodriguez:Sottile:2018}. It can be shown that no proper subset of the points in $W \cap M_t$ is an affine linear function of $t$. 

This leads to the idea of the trace test: Let $t_1 \in \mathbb{C}\setminus \{0\}$, fix $W' \subseteq W \cap M_0$ and compute
${\rm tr}(t_1) := (w(t_1) - w(0)) - (w(0) - w(-t_1))$.
Note that ${\rm tr}(t_1)$ is identically zero if and only if $w$ is an affine linear function of $t$, which is true if and only if the cardinality of $W'$ corresponds to the degree of $W$.
Due to the generality assumption on $M_t$ it is sufficient to compute ${\rm tr}(t_1)$ for only \emph{one} $t_1 \in \mathbb{C}\setminus \{0\}$.

\section{A Numerical Approach}\label{sec:results}

In this section we explain our use of numerical algebraic geometry to obtain Theorem* \ref{thr:main} below.
Reasonable mathematicians may differ as to whether it is appropriate to state this result as a theorem since we currently cannot certify the last step of our computation.
We add the asterisk to acknowledge these differing opinions. 

\begin{theorem*}\label{thr:main}
  The degree of the orbit closure of a general cubic surface under the action of \(\pgl(\comp,4)\) is 96120.
\end{theorem*}
\noindent All computations performed to arrive at this result are available from the authors upon request. 

To compute the degree of the orbit closure, we sample a general cubic surface $f \in \PC^{19}$ by drawing the real and imaginary parts of each of its coordinates independently from a univariate normal distribution. 
We then solve the polynomial system \eqref{eq:formulation_enumerative} encoding the enumerative geometry problem.
A naive strategy is to sample 15 points $p_1,\ldots,p_{15} \in \PC^{3}$ in general position and use a total degree homotopy, but in this case the B\'ezout bound is $3^{15}=14,348,907$.
Here, the monodromy method is substantially more efficient.

To apply the monodromy method, we consider \eqref{eq:formulation_enumerative} as a polynomial system on the entries of \(\varphi\) parameterized by 15 points $p_1,\ldots,p_{15}$ in $\PC^3$.
We consider the incidence variety
\begin{equation*}\label{eq:solution-variety}
  V = \{ (\varphi, (p_1,\ldots,p_{15}))\in \PC^{15}\times (\PC^3)^{15} \; | \; F(\varphi(p_i)) = 0 , \, i=1,\ldots,15  \} \,
\end{equation*}
and we denote by $\pi$ the projection \(\PC^{15}\times(\PC^{3})^{15}\to (\PC^3)^{15}\) restricted to $V$.

We find a start pair $(\varphi_0; p_1,\dots,p_{15}) \in V$ and then we use the monodromy action on the fiber $\pi^{-1}(p_1,\dots,p_{15})$ to find all solutions in this fiber.
Such a start pair can be found by exchanging the role of variables and parameters. First, we sample a $\varphi_0 \in \PC^{15}$ and the first three coordinates of 15 points $p_i \in \PC^3$ in general position.
This yields a system of 15 polynomials each depending only on a unique variable: The $i$th polynomial depends only on the fourth coordinate of  $p_i$.
Such a system is easy to solve.
Solving it yields a start pair $(\varphi_0; p_1,\dots,p_{15}) \in V$, on which we run the monodromy method implemented in the software package {\tt HomotopyContinuation.jl} \cite{HomotopyContinuation.jl}.
In less than an hour on a single core, this method found 96120 approximate solutions corresponding to the start points $p_1,\ldots, p_{15} \in \PC^3$.

Next we apply Smale's \(\alpha\)-theory as implemented in the software {\tt alphaCertified} \cite{Hauenstein:Sottile:2012} to certify two conditions of our numerical approximations: First, we show that each is indeed an approximate zero to our original polynomial system, and second that all 96120 approximate zeros have distinct associated zeros.
Due to computational limits we were only able to obtain a certificate using (arbitrary precision) floating point arithmetic.
Hauenstein and Sottile call this a ``soft'' certificate since it does not eliminate the possibility of floating point errors.
It is preferable to use rational arithmetic for certification, but for a system of our size too much time is required to perform such a computation.

The certification process establishes a \emph{lower} bound on the degree of the orbit closure. 
As a last step, we run a trace test to verify that we have indeed found \emph{all} solutions.
The trace test described in the previous section is only applicable to varieties $W \subset \PC^n$. In \cite{Leykin:Rodriguez:Sottile:2018} the authors derive a trace test to certify the completeness of a \emph{collection of partial multihomogeneous witness sets}. Our formulation \eqref{eq:formulation_enumerative} provides only one partial multihomogeneous witness set, namely $\pi^{-1}(p_1,\dots,p_{15})$,  and not the entire collection that would be necessary to run a trace test. 
To avoid these complications, we use formulation \eqref{eq:formulation_subspace}.
We note that it is straightforward to construct a linear subspace $L$ from the 15 points $p_1,\ldots, p_{15}$ such that our solutions from the monodromy computation are also solutions to \eqref{eq:formulation_subspace}, so translating formulation  \eqref{eq:formulation_enumerative} to \eqref{eq:formulation_subspace} is not difficult.

In the language of numerical algebraic geometry our 96120 solutions together with the linear subspace $L$ constitute a \emph{pseudo witness set} \cite{Hauenstein:Sommese:2010}.
We construct a general pencil $M_t$ of linear spaces with $M_0 = L$.
Working with approximate solutions refined to around 38 digits of accuracy we obtain for ${\rm tr}(1)$ a vector with norm of approximately $10^{-32}$.
Additionally, increasing the accuracy of the solutions decreases the norm of the trace test result.
While this gives us very high certainty that we indeed obtained all solutions, we do not have a rigorous certificate that the trace test converges to zero when we increase the accuracy of the solutions.
A certification of the trace test similar to Smale's $\alpha$-theory for numerical solutions remains an open problem. 

From the described computations we conclude that degree of the orbit closure of a general cubic surface under the action \(\pgl(\comp,4)\) is 96120.

We note that as a test of our methods, we confirmed known degrees of other varieties.
In agreement with a theoretical result of Aluffi and Faber  \cite{Aluffi:Faber:1993}, we  computed that the degree of the orbit closure of a general quartic curve in the plane is 14280.
Additionally we computed that the degree of the orbit closure of the Cayley cubic, defined by the equation $yzw+xzw+xyw+xyz=0$, is $305$.
Due to the symmetry of the variables in the Cayley cubic, there are $4!$ matrices corresponding to every polynomial in the orbit.
As expected, we computed $7320 = 4! \cdot 305$ solutions.
This coincides with a theoretical result of Vainsencher \cite{cayleycubic}.

\setcounter{tocdepth}{1}

\bibliography{sample}

\end{document}